\documentclass[12pt]{article}

\newtheorem{definition}{{\bf Definition}}[section]
\newtheorem{theorem}[definition]{{\bf Theorem}}

\newtheorem{proposition}[definition]{{\bf Proposition}}

\def\ffi{\varphi}
\def\<{\langle}
\def\>{\rangle}

\def\ot{\otimes}

\def\tr{{\rm tr}}
\def\Mn{{M_n({\mathbb C})}}
\def\M2{{M_n({\mathbb C}) \otimes M_n({\mathbb C})}}
\def\rank{{\rm rank}}
\def\dis{\displaystyle}

\usepackage{amscd,amsmath,amssymb,graphicx}

\begin{document}

\begin{center}
 {\Large Maximal rank of extremal marginal tracial states}@
\end{center}

\begin{center}
Hiromichi Ohno \\
Department of Mathematics, Faculty of Engineering, Shinshu University,\\
4-17-1 Wakasato, Nagano 380-8553, Japan \\
\end{center}

\section{Introduction}
States on coupled quantum systems are studied in many point of view.
For example, entangled states and separable states are
investigated in many papers. 
In this paper, we consider states on coupled quantum systems $\M2$ whose
restrictions to each subsystems are the normalized traces.
Such states are called marginal tracial states.
In \cite{Pa1}, Parthasarathy showed that
every extremal marginal tracial state on 
$M_2({\mathbb C})\ot M_2({\mathbb C})$ is a pure state.
In \cite{PS}, Price and Sakai introduced
necessary and sufficient conditions for a state to be extremal
marginal tracial states.
This problem is an analogue of the Birkhoff Theorem which says
extremal points of a set of doubly stochastic matrices of order $n$ are
the permutation matrices of order $n$.

Price and Sakai conjectured that 
every extremal marginal tracial state is pure if $n\ge 3$ in \cite{PS}.
But by using the one-to-one correspondence between 
marginal tracial states on $\M2$ and unital completely positive trace
preserving (UCPT) maps on $\Mn$ (see e.g. \cite{Ar2,Ru}), and the result of Landau 
and Streater \cite[Theorem 1]{LS} which shows there exist nonunitary
extremal UCPT maps if $n \ge 3$, we can conclude that 
there exist nonpure extremal marginal tracial states.
This result, first known by Arveson, was told by Price in
private communication.
Since there is not a paper which is written about this result,
we see this fact in this paper for completeness.
Moreover, we show the maximal rank of extremal marginal tracial states for
some special cases and consider diagonal marginal tracial states which
correspond to diagonal UCPT maps introduced in \cite{LS}.

In Section \ref{sect2}, we introduce the relation between 
marginal taracial states and UCPT maps and 
construct nonpure extremal marginal tracial
states.
Furthermore, we fix the maximal rank of extremal marginal tracial states on
$\M2$ when $n=3,4$ and see a lower bound of maximal rank of extremal 
marginal tracial states for $n \ge 5$.
In Section \ref{sect3}, the maximal rank of extremal diagonal marginal tracial
states are investigated. Moreover, 
we show that extremal diagonal marginal tracial states are dense in a set of 
all diagonal marginal tracial states.


\section{Marginal tracial states and UCPT maps}\label{sect2}

In this section, we consider a relation between 
marginal tracial states and UCPT (unital completely positive trace
preserving) maps, and the maximal rank of extremal marginal tracial states.

\begin{definition}
A state $\rho$ on $\M2$ is a marginal tracial state 
if and only if the restrictions of $\rho$ to
$\Mn \otimes I$ and $I \otimes \Mn$ are the normalized traces.
$\Gamma(n)$ is a set of all marginal tracial states on $\M2$.
\end{definition}

Extremal marginal tracial states are considered in \cite{Pa1, PS}
and
Parthasarathy proved that an extremal marginal tracial state is 
a pure state, if $n=2$ (\cite{Pa1}).
We will show that this is not true when $n \ge 3$ by using a one-to-one
correspondence  between marginal tracial states and UCPT maps.

We denote $UCPT(n)$ is a set of all UCPT maps from
$\Mn$ to $\Mn$. A UCPT map $\ffi \in UCPT(n)$ can be written as
\[
\ffi (A) = \sum_{i=1}^k v_i^* A v_i
\]
for some $\{v_i\} \subset \Mn$ with 
$\sum_{i=1}^k v_i^*v_i = \sum_{i=1}^k v_i v_i^* =I$.
This expression is not unique. But the matrices $\{v_i\}$ can be taken
linearly independent, and then the number $k$  of terms is uniquely determined (see \cite{Choi})
and we denote $r(\ffi) =k$.

For any $\ffi \in UCPT(n)$, we can define a state $\pi(\ffi)$
on $\M2$ by
\[
\pi(\ffi)(A\ot B) = \tr(\ffi(A) {}^tB) = \< (\ffi(A) \ot B) \xi , \xi \>,
\]
where $\tr$ is the normalized trace, ${}^tB$ is a transpose of $B$ and
$\xi = {1\over \sqrt{n}}
\sum_{i=1}^n e_i \ot e_i \in {\mathbb C}^n \ot {\mathbb C}^n$
with a fixed orthonormal basis $\{e_i\}$ of ${\mathbb C}^n$.
Since there is a one-to-one correspondence between states and density matrices,
we will write $D_{\rho}$ as the corresponding 
density matrix of a state $\rho$ and write ${\rm rank}(\rho) := 
{\rm rank}(D_{\rho})$.
Here 
\begin{equation}\label{density}
D_{\pi(\ffi)} = \sum_{i=1}^k |(v_i \ot I) \xi \>\< (v_i \ot I)\xi| 
={1\over n} \sum_{i,j=1}^n \ffi(e_{ij})\ot e_{ij},
\end{equation}
where $\{e_{ij}\}$ is a set of matrix units of $\Mn$,
so that this is a Choi matrix.
This correspondence is well known but we 
show the proof of next theorem for completeness.

\begin{theorem}
A map $\pi$ from $UCPT(n)$ to $\Gamma(n)$ defined above is bijective.
In particular, $\pi(\ffi)$ is an extremal point of $\Gamma(n)$ if and only if 
$\ffi$ is an extremal point of $UCPT(n)$.
Moreover, ${\rm rank}(\pi(\ffi))$ is equal to $r(\ffi)$.
\end{theorem}
\begin{proof}
First, $\pi(\ffi)$ is a marginal tracial state for any $\ffi \in 
UCPT(n)$. Indeed, 
\begin{eqnarray*}
\pi(\ffi) (A\ot I ) &=& \tr(\ffi(A) ^{t}I) = \tr (\ffi(A)) = \tr(A), \\
\pi(\ffi) (I\ot B ) &=& \tr(\ffi(I) ^tB) = \tr( ^tB) = \tr(B).
\end{eqnarray*}

Next, we show that $\pi$ is surjective.
For a marginal tracial state $\rho$, the spectral decomposition of the
density matrix have a form $D_\rho = \sum_{i=1}^k \lambda_i |\zeta_i\>\<
\zeta_i|$ for some $\lambda_i \ge0$ and 
$\zeta_i \in {\mathbb C}^{n^2} = {\mathbb C}^n \ot {\mathbb C}^n$.
We write
\[
\zeta_i = \sum_{j=1}^n \zeta_{ij}\ot e_j \in {\mathbb C}^n \ot {\mathbb C}^n.
\]
Define $v_i \in \Mn$ by $v_i(e_j) = \sqrt{n} \zeta_{ij}$.
Let $\ffi = \sum_{i=1}^k {\lambda_i} 
v_i^* \, \cdot \, v_i$, then $\pi(\ffi) = \rho$ and $\ffi \in UCPT(n)$.
Indeed, since we can extend the domain of $\pi$ to 
all completely positive maps and $\ffi$ is a completely positive map,
we can define a positive linear functional $\pi (\ffi)$ on $\M2$.
Then by (\ref{density}),
\begin{eqnarray*}
D_{\pi(\ffi)} &=&
 \sum_{i=1}^k \lambda_i |(v_i \ot I) \xi \>\< (v_i\ot I)\xi|
=\sum_{i=1}^k \lambda_i \big| \sum_{j=1}^k (\zeta_{ij}\ot e_j) \big\>
\big\< 
\sum_{j=1}^k (\zeta_{ij}\ot e_j)\big| \\
&=& \sum_{i=1}^k \lambda_i |\zeta_i \>\< \zeta_i| 
= D_\rho.
\end{eqnarray*}
Moreover we have
\begin{eqnarray*}
\tr(\ffi(A)) &=& \tr(\ffi(A) {}^tI) = \rho(A \otimes I ) = \tr (A),\\
\tr(\ffi(I) {}^tB) &=& \rho(I\otimes B) = \tr(B) = \tr ( ^tB)
\end{eqnarray*}
for all $A, B\in \Mn$. This implies that
$\ffi$ is a UCPT map so that $\pi$ is surjective.
It is easy to see that $\pi$ is injective by definition.
Therefore $\pi$ is bijective.

Finally, we see ${\rm rank}(\pi(\ffi)) = r(\ffi)$.
Since $\{ v_i \}$ is linearly independent, $\{(v_i \ot I) \xi \}$ is also
linearly independent. Hence the density matrix
\[
D_{\pi(\ffi)} = \sum_{i=1}^k |(v_i \ot I )\xi\>\< (v_i \ot I )\xi|
\]
is rank $k$. This implies ${\rm rank}(\pi(\ffi)) = r(\ffi)$.
\end{proof}

By Theorem 1 in \cite{LS}, 
there exists an extremal point $\ffi$ in $UCPT(n)$ 
with $r(\ffi) \ge 2$ (so $\ffi$ is a nonunitary map), if $n \ge 3$.
This implies that there exists an extremal marginal 
tracial state in $\Gamma(n)$ which is not a pure state, if $n \ge 3$.
Moreover, we will use the following theorem to construct examples.

\begin{theorem}[\cite{LS}]\label{extremalucpt}
Let $\ffi = \sum_{i=1}^k v_i^* \, \cdot \, v_i$ and $\{ v_i \}$ is linearly 
independent with $\sum_{i=1}^k v_i^* v_i = \sum_{i=1}^k v_i v_i^* = I$. 
Then $\ffi$ is an extremal point of $UCPT(n)$ if and only if
$\{v_i v_j^*\}_{i,j=1}^k$ and $\{v_j^* v_i\}_{i,j=1}^k$
are bi-independent sets.
\end{theorem}

Next we consider the maximal rank of extremal marginal tracial states
$\Gamma(n)$ which is denoted by $MR(n)$.
From the next theorem proven in \cite{PS}, 
we can obtain the upper bound of 
$MR(n)$. 

\begin{theorem}[\cite{PS}]
Let $\rho$ be a marginal tracial state on $\M2$ and $P_\rho$ be a
support projection of $\rho$. Then the following conditions are equivalent.

(i) $\rho$ is an extremal point of $\Gamma(n)$,

(ii) $(P_\rho (\M2)P_\rho) \cap ((\Mn \ominus {\mathbb C}I) 
\otimes (\Mn \ominus {\mathbb C}I)) = \{ 0\}$.
\end{theorem}

Since the dimensions of $P_\rho (\M2)P_\rho$, $(\Mn \ominus {\mathbb C}I) 
\otimes (\Mn \ominus {\mathbb C}I)$ and $\M2$ are 
$\rank(\rho)^2$, $(n^2-1)^2$ and $n^4$, respectively, we obtain 
\begin{equation}\label{upperbound}
MR(n) \le \sqrt{2n^2-1}.
\end{equation}

We know that extremal points of $\Gamma(2)$ are pure so that $MR(2)=1$.
The following theorems fix $MR(3)$ and $MR(4)$.

\begin{theorem}
The maximal rank of extremal marginal tracial states $\Gamma(3)$ is $4$.
\end{theorem}
\begin{proof}
From (\ref{upperbound}), $MR(3) \le 4$. Hence we only need to
construct a UCPT map $\ffi$ which is extremal in $UCPT(3)$ and $r(\ffi)=4$.
Let 
\begin{eqnarray*}
w_1 = e_{11}, \quad w_2 = e_{12}+\sqrt{2} e_{23}, \quad
w_3 =\sqrt{2}e_{21} + \sqrt{3} e_{32},\quad
w_4 =e_{31} + \sqrt{2}e_{13}.
\end{eqnarray*}
Then we can see that 
\[
\sum_{i=1}^4 w_i w_i^* =\sum_{i=1}^4 w_i^* w_i = 4 I.
\]
Hence $v_i = {1\over 2} w_i$ satisfies $\sum_{i=1}^4 
v_i v_i^* = \sum_{i=1}^4 v_i^* v_i = I $.
Moreover $\{w_i w_j^*\}$ and $\{w_j^* w_i\}$ are bi-independent sets. 
Indeed,
since we have
\begin{eqnarray*}
\begin{array}{llll}
w_1w_1^* = e_{11} & w_2w_1^* = 0 & 
w_3 w_1^* =\sqrt{2}e_{21} & w_4 w_1^* = e_{31} \\
w_1w_2^* = 0 & w_2w_2^* = e_{11}+2e_{22} & 
w_3 w_2^* =\sqrt{3}e_{31} & w_4 w_2^* = 2e_{12} \\
w_1w_3^* = \sqrt{2}e_{12} & w_2w_3^* = \sqrt{3}e_{13} & 
w_3 w_3^* =2e_{22} +3e_{33} & w_4 w_3^* = \sqrt{2}e_{32} \\
w_1w_4^* = e_{13} & w_2w_4^* = 2e_{21} & 
w_3 w_4^* =\sqrt{2}e_{23} & w_4 w_4^* = e_{33}+2e_{11} 
\end{array}
\end{eqnarray*}
and
\begin{eqnarray*}
\begin{array}{llll}
w_1^*w_1 = e_{11} & w_2^*w_1 = e_{21} & 
w_3^* w_1 =\sqrt{2}e_{12} & w_4^* w_1 = \sqrt{2}e_{31} \\
w_1^*w_2 = e_{12} & w_2^*w_2 = e_{22}+2e_{33} & 
w_3^* w_2 =2e_{13} & w_4^* w_2 = \sqrt{2}e_{32} \\
w_1^*w_3 = 0 & w_2^*w_3 = 2e_{31} & 
w_3^* w_3 =2e_{11} +3e_{22} & w_4^* w_3 = \sqrt{3}e_{12} \\
w_1^*w_4 = \sqrt{2}e_{13} & w_2^*w_4 = \sqrt{2}e_{23} & 
w_3^* w_4 =\sqrt{3}e_{21} & w_4^* w_4 = e_{11}+2e_{33}, 
\end{array}
\end{eqnarray*}
$\sum_{i,j=1}^4 a_{ij} v_iv_j^* = \sum_{i,j=1}^4 a_{ij} v_j^* v_i =0$
imply 
$a_{ij} = 0$ for all $1\le i,j \le 4$ from a simple calculation so that
$\{w_i w_j^*\}$ and $\{w_j^* w_i\}$ are bi-independent sets.
Therefore $\{v_i v_j^*\}$ and $\{v_j^* v_i\}$ are bi-independent sets
and $\ffi$ is an extremal point of $UCPT(3)$ with $r(\ffi) =4$ by theorem
\ref{extremalucpt}.
This shows $MR(3) = 4$.
\end{proof}

\begin{theorem}
The maximal rank of extremal marginal tracial states $\Gamma(4)$ is $5$.
\end{theorem}
\begin{proof}
From (\ref{upperbound}), $MR(4) \le 5$. Hence we only need to
construct $\ffi$ which is extremal in $UCPT(4)$ and $r(\ffi)=5$.
Let 
\begin{eqnarray*}
w_1 &=& e_{13} + e_{32}, \quad w_2 = \sqrt{2}e_{24}+ \sqrt{2}e_{43}, \quad
w_3 =\sqrt{2}e_{14} + \sqrt{3} e_{31},\\
w_4 &=& e_{21} + \sqrt{2}e_{42}, \quad w_5 = e_{12}+e_{23}.
\end{eqnarray*}
Then we can see that 
\[
\sum_{i=1}^5 w_i w_i^* = \sum_{i=1}^5 w_i^* w_i = 4 I.
\]
Hence $v_i = {1\over 2} w_i$ satisfy $\sum_{i=1}^5 
v_i v_i^* = \sum_{i=1}^5 v_i^* v_i = I $.
Moreover $\{w_i w_j^*\}$ and $\{w_j^* w_i\}$ are bi-independent sets.
 Indeed,
since we have
\begin{eqnarray*}
&&\begin{array}{lll}
w_1w_1^* = e_{11}+e_{33} & 
w_2w_1^* = \sqrt{2}e_{41} & 
w_3 w_1^* = 0  \\
w_1w_2^* = \sqrt{2}e_{14} & 
w_2w_2^* = 2e_{22}+2e_{44} & 
w_3 w_2^* =2e_{12} \\
w_1w_3^* = 0 & 
w_2w_3^* = 2e_{21} & 
w_3 w_3^* =2e_{11} +3e_{33}  \\
w_1w_4^* = \sqrt{2}e_{34} & 
w_2w_4^* = 0 & 
w_3 w_4^* =\sqrt{3}e_{32}\\
w_1w_5^* = e_{31}+e_{12} & 
w_2w_5^* = \sqrt{2}e_{42} & 
w_3 w_5^* = 0
\end{array} \\
&&\begin{array}{ll}
w_4 w_1^* = \sqrt{2} e_{43} &
w_5 w_1^* = e_{13}+e_{21} \\
w_4 w_2^* = 0 & 
w_5w_2^* = \sqrt{2}e_{24} \\
w_4 w_3^* = \sqrt{3}e_{23} &
w_5w_3^* = 0 \\
w_4 w_4^* = e_{22}+2e_{44} &
w_5w_4^* = \sqrt{2} e_{14} \\ 
w_4 w_5^* = \sqrt{2} e_{41} &
w_5w_5^* = e_{11} + e_{22} \\ 
\end{array}
\end{eqnarray*}
and
\begin{eqnarray*}
&&
\begin{array}{lllll}
w_1^*w_1 = e_{33} +e_{22}&
w_2^*w_1 = 0 & 
w_3^* w_1 =\sqrt{2}e_{43} + \sqrt{3}e_{12} \\
w_1^*w_2 = 0 & 
w_2^*w_2 = 2e_{44}+2e_{33} & 
w_3^* w_2 = 0  \\
w_1^*w_3 = \sqrt{2}e_{34} + \sqrt{3}e_{21} & 
w_2^*w_3 = 0 & 
w_3^* w_3 =2e_{44} +3e_{11}  \\
w_1^*w_4 = 0 & 
w_2^*w_4 = \sqrt{2}e_{41} + 2e_{32} & 
w_3^* w_4 = 0  \\ 
w_1^*w_5 = e_{32} & 
w_2^*w_5 = \sqrt{2}e_{43} & 
w_3^* w_5 = \sqrt{2}e_{42}
\end{array}
\\
&&
\begin{array}{ll}
w_4^* w_1 = 0 &
w_5^*w_1 = e_{23} \\
w_4^* w_2 = \sqrt{2}e_{14} + 2e_{23} &
w_5^*w_2 = \sqrt{2}e_{34} \\
w_4^* w_3 = 0 &
w_5^* w_3 = \sqrt{2} e_{24} \\
w_4^* w_4 = e_{11}+2e_{22} &
w_5^* w_4 = e_{31} \\ 
w_4^* w_5 = e_{13} &
w_5^* w_5 = e_{22}+e_{33}, 
\end{array}
\end{eqnarray*}
$\sum_{i,j=1}^5 a_{ij} v_iv_j^* = \sum_{i,j=1}^5 a_{ij} v_j^* v_i =0$
imply 
$a_{ij} = 0$ for all $1\le i,j \le 5$ from a simple calculation so that
$\{w_i w_j^*\}$ and $\{w_j^* w_i\}$ are bi-independent sets.
Therefore $\{v_i v_j^*\}$ and $\{v_j^* v_i\}$ are bi-independent sets
and $\ffi$ is an extremal point of $UCPT(4)$ with $r(\ffi) =5$ by theorem
\ref{extremalucpt}.
This shows $MR(4) = 5$.
\end{proof}

The next theorem shows a lower bound of $MR(n)$ for $n \ge 5$

\begin{theorem}
The maximal rank of extremal marginal tracial states $\Gamma(n)$ is 
at least $n$.
\end{theorem}
\begin{proof}
We construct $\ffi$ which is extremal in $UCPT(n)$ and 
$r(\ffi)=n$ for $n\ge 5$.
Let 
\begin{eqnarray*}
v_1  &=& \sqrt{n-2\over n-1} \sum_{j=2}^n e_{jj}  \\
v_i  &=& {1\over \sqrt{n-1}} \left(e_{1i} + e_{i1}\right)
\end{eqnarray*}
for $2 \le i \le n$. Then we have 
\[
\begin{array}{ll}
v_1^*v_1 = \dis{n-2\over n-1} \sum_{j=2}^n e_{jj}, &  \\
v_1^* v_j = \dis{\sqrt{n-2}\over n-1}e_{j1}, & {\rm for} \,\, j\ge 2, \\
v_j^* v_1 = \dis{\sqrt{n-2}\over n-1}e_{1j}, & {\rm for} \,\, j\ge 2, \\
v_j^* v_j = \dis{1\over n-1}e_{11} + e_{jj}, &  {\rm for}\,\, j\ge 2, \\
v_j^* v_k = \dis{1\over n-1}e_{jk}, & 
{\rm for}\,\, j,k\ge 2 \,\,{\rm and} \,\, j\neq k.
\end{array}
\]
From a simple calculation, we obtain that $\{v_j^* v_i\}$
 is linearly independent
and so that 
$\ffi$ is an extremal point of $UCPT(n)$ with $r(\ffi) =n$ by theorem
\ref{extremalucpt}.
This shows $MR(n) \ge n$ for $n\ge 5$.
\end{proof}


\section{diagonal UCPT maps}\label{sect3}

In this section, we consider diagonal UCPT maps.
Diagonal maps are introduced in \cite{LS}.

\begin{definition}[\cite{LS}]
A linear map $\ffi$ from $\Mn$ to $\Mn$ is diagonal
if it has a form
\[
\ffi(A) = C \circ A
\]
for some $C \in \Mn$, where $C \circ A$ is the Schur product of $C$ and $A$.
\end{definition}

A marginal tracial state corresponding to a diagonal UCPT map
is called a diagonal marginal tracial state.

Completely positive diagonal maps are characterized by the next proposition.

\begin{proposition}[\cite{LS}]
$\ffi$ is a completely positive diagonal map, 
if and only if in any representation
\[
\ffi = \sum_{i=1}^k v_i^* \,\cdot \, v_i,
\]
the matrices $v_i$ are diagonal.
\end{proposition}

In \cite{LS}, it is also shown that the maximal rank of 
diagonal UCPT maps on $\Mn$ is at most $\sqrt{n}$.
The next theorem shows that for any $a^2 \le n$ 
we can construct a diagonal UCPT map $\ffi$ with
$r(\ffi) = a$.

\begin{theorem}\label{thm3.3}
Let $a \in {\mathbb N}$ be such that $a^2 \le n$ for some
$n\ge 4$. Then there are diagonal matrices 
$v_1, \ldots, v_a \in \Mn$ such that the map
\[
\ffi(A) = \sum_{i=1}^a v_i^* A v_i, \qquad A \in \Mn
\]
is an extremal UCPT map, hence $\ffi$ corresponds to an
extremal diagonal marginal tracial state on $\M2$ of rank $a$.
\end{theorem}

\begin{proof}
The result is trivial if $a=1$: let $v_1$ be any diagonal unitary matrix 
in $\Mn$. So we shall suppose $a \ge 2$.
For notational purposes let us agree to write our diagonal matrices $v$ as
row vector with $n$ entries.

Consider the following vectors of length $m=a^2$:
\begin{eqnarray*}
v_1 &=& (\overbrace{1,0,0, \ldots, 0,0}^{{\rm length} \,\, a}, 
b_{11}, b_{12},\ldots
b_{1l}), \\
v_2 &=& ({0,1,0 \ldots, 0,0}, 
b_{21}, b_{22},\ldots
b_{2l}), \\
&\vdots& \\
v_a &=& ({0,0,0 \ldots, 0,1}, 
b_{a1}, b_{a2},\ldots
b_{al}),
\end{eqnarray*}
where $l = m-a = a^2 -a$.
We will construct entries $b_{ij}$, $1\le i,j \le l$, such that
the set $\{v_iv_j^*\}_{i,j=1}^a$ is linearly independent and such that
$\sum_{i=1}^a v_i^* v_i = \sum_{i=1}^a v_iv_i^*=I$.
This will make the completely positive map
$\ffi = \sum_{i=1}^a v_i^* \, \cdot \, v_i$ 
a UCPT map which is extremal in $UCPT(m)$.

Now choose
$\theta_1, \ldots, \theta_a \in [0,2\pi]$ such that the set of 
elements $\{\theta_i - \theta_j \,:\, 1 \le i \neq j \le a \}$
are all distinct mod $2\pi$.
Let 
\begin{eqnarray*}
\begin{array}{cccc}
b_{11} = \dis{1\over \sqrt{a}} e^{i\theta_1} &
b_{12} = \dis{1\over \sqrt{a}} e^{2i\theta_1} &
\ldots &
b_{1l} = \dis{1\over \sqrt{a}} e^{li\theta_1}  \\
b_{21} = \dis{1\over \sqrt{a}} e^{i\theta_2} &
b_{22} = \dis{1\over \sqrt{a}} e^{2i\theta_2} &
\ldots &
b_{2l} = \dis{1\over \sqrt{a}} e^{li\theta_2}  \\
\vdots & \vdots & \vdots & \vdots \\
b_{a1} = \dis{1\over \sqrt{a}} e^{i\theta_a} &
b_{a2} = \dis{1\over \sqrt{a}} e^{2i\theta_a} &
\ldots &
b_{al} = \dis{1\over \sqrt{a}} e^{li\theta_a}.  \\
\end{array}
\end{eqnarray*}
Write the set $\{v_i v_j^*\}_{1\le i,j \le a}$ as a list of 
$a^2$ row vectors, $v_1v_1^*, v_1v_2^*, \ldots, v_1v_a^*, 
v_2v_1^*, v_2v_2^*, \ldots v_2v_a^*, \ldots, v_av_a^*$, which looks like
\begin{eqnarray*}
v_1v_1^* &=& (\overbrace{1,0, 0,  \ldots, 0,0}^{{\rm length}\,\, a}, 
\overbrace{\dis{1\over a}, \dis{1\over a}, \ldots, 
\dis{1\over a}}^{{\rm length}\,\,l} )\\
v_1v_2^* &=& (0,0, 0,  \ldots, 0,0, 
{1\over a}e^{i(\theta_1 -\theta_2)}, {1\over a}e^{2i(\theta_1 -\theta_2)}, 
\ldots, {1\over a}e^{li(\theta_1 -\theta_2)} )\\
v_1v_3^* &=& (0,0, 0,  \ldots, 0,0, 
{1\over a}e^{i(\theta_1 -\theta_3)}, {1\over a}e^{2i(\theta_1 -\theta_3)}, 
\ldots, {1\over a}e^{li(\theta_1 -\theta_3)} )\\
&& \vdots \\
v_av_a^* &=& (0,0, 0,  \ldots, 0,1, 
{1\over a}, {1\over a}, 
\ldots, {1\over a} ).
\end{eqnarray*}
Rearrange the rows so that $v_1v_1^*, v_2v_2^* \ldots v_av_a^*$ 
appear first
and one obtains an $m \times m$ matrix that  looks like
\[
\left[
\begin{array}{ccccccccccc}
1 & 0 & 0 & \cdots & 0 & 0 &\dis{1\over a} &\dis{1\over a} &\dis{1\over a} 
& \cdots &\dis{1\over a} \\
0 & 1 & 0 & \cdots & 0 & 0 &\dis{1\over a} &\dis{1\over a} &\dis{1\over a} 
& \cdots &\dis{1\over a} \\
\vdots & \vdots & \vdots & \vdots & \vdots & \vdots & \vdots & \vdots & 
\vdots & \vdots & \vdots \\
0 & 0 & 0 & \cdots & 0 & 1 &\dis{1\over a} &\dis{1\over a} &\dis{1\over a} 
& \cdots &\dis{1\over a} \\
0 & 0 & 0 & \cdots & 0 & 0 &
\dis{1\over a}e^{i(\theta_1 -\theta_2)} &
\dis{1\over a}e^{2i(\theta_1 -\theta_2)} &
\dis{1\over a}e^{3i(\theta_1 -\theta_2)} 
& \cdots & \dis{1\over a}e^{li(\theta_1 -\theta_2)} \\
\vdots & \vdots & \vdots & \vdots & \vdots & \vdots & \vdots & \vdots & 
\vdots & \vdots & \vdots \\
0 & 0 & 0 & \cdots & 0 & 0 &
\dis{1\over a}e^{i(\theta_{a} -\theta_{a-1})} &
\dis{1\over a}e^{2i(\theta_a -\theta_{a-1})} &
\dis{1\over a}e^{3i(\theta_a -\theta_{a-1})} 
& \cdots & \dis{1\over a}e^{li(\theta_a -\theta_{a-1})} 
\end{array}
\right].
\]
Note that this matrix has rank $m$ if and only if 
the lower corner $l \times l$ matrix has rank $l$.
But if one factors out $1/a$ from each row of this 
$l \times l$ matrix then the remaining matrix looks like
\[
\left[
\begin{array}{cccc}
e^{i(\theta_{1} -\theta_{2})} & e^{2i(\theta_{1} -\theta_{2})} &
\cdots &e^{li(\theta_{1} -\theta_{2})} \\
e^{i(\theta_{1} -\theta_{3})} & e^{2i(\theta_{1} -\theta_{3})} &
\cdots &e^{li(\theta_{1} -\theta_{3})} \\
\vdots & \vdots & \vdots & \vdots \\
e^{i(\theta_{a} -\theta_{a-1})} & e^{2i(\theta_{a} -\theta_{a-1})} &
\cdots &e^{li(\theta_{a} -\theta_{a-1})}
\end{array}
\right].
\]
Note that this matrix is a Vandermonde matrix with non-zero 
determinant because of the way that the numbers $\theta_1, \ldots,
\theta_a$ were chosen.
This shows that the set $\{v_iv_j^*\}_{1\le i, j \le a}$ is a
linearly independent set.
Since $\sum_{i=1}^a v_i^* v_i = \sum_{i=1}^a v_i v_i^* = I$,
we are done if $a = n^2$.

If $a^2 < n$, then add $n-a^2$ entries to the end of each vectors 
$v_1, \ldots , v_a$ to form
\begin{eqnarray*}
v_1' &=& v_1 \oplus (c_{11} , c_{12}, \ldots , c_{1, n-m}) \\
v_2' &=& v_2 \oplus (c_{21} , c_{22}, \ldots , c_{2, n-m}) \\
&\vdots& \\
v_a' &=& v_a \oplus (c_{a1} , c_{a2}, \ldots , c_{a, n-m}) 
\end{eqnarray*}
in such a way that for all $1\le j \le n-m$,
$\sum_{i=1}^a |c_{ij}|^2 = I$.
This will give us 
\[
\sum_{i=1}^a {v_i'}^{*} v_i' = \sum_{i=1}^a v_i' {v_i'}^{*} =I.
\]
Furthermore, since $\{v_iv_j^*\}_{1\le i,j \le a}$ is
a linearly independent set, then so is 
 $\{v_i'{v_j'}^{*}\}_{1\le i,j \le a}$, and therefore the mapping
\[
\ffi' (A) = \sum_{i=1}^a {v_i'}^{*} A v_i', \qquad A \in \Mn
\]
is an extremal UCPT map, hence corresponds to an extremal 
marginal diagonal tracial state on $\Mn$ of rank $a$.
\end{proof}

Moreover, we can show that extremal diagonal UCPT maps are dense
in the set of all diagonal UCPT maps.

\begin{theorem}
Extremal diagonal marginal tracial states of rank $a$ are dense 
inside the set of all diagonal marginal tracial states of rank $a$ or less.
\end{theorem}

\begin{proof}
Let $\ffi$ be diagonal UCPT maps on $M_n({\mathbb C})$ with
\[
\ffi(A) = \sum_{i=1}^a u_i^* A u_i, \qquad 
\]
for any $A \in M_n({\mathbb C})$ which is allowed that $u_i = 0$,
and let $\psi$ be diagonal UCPT maps on $M_n({\mathbb C})$ with
\[
\psi(A) = \sum_{i=1}^a v_i^* A v_i
\]
for any $A \in M_n({\mathbb C})$.
Moreover we assume that $\{v_i v_j^*\}$ is linearly independent so that
$\psi$ is an extremal point in the set of all UCPT maps.
Let $w_i = u_i + \varepsilon v_i$ and 
\[
\ffi_\varepsilon (A) :=  (\sum_{i =1}^a w_i^* w_i)^{-1} \sum_{i=1}^a w_i^* A w_i
\]
for all $A \in M_n({\mathbb C})$. Since $\sum_{i =1}^a u_i^* u_i = I$,
$\sum_{i =1}^a w_i^* w_i$ is invertible for small $\varepsilon$.
Therefore $\ffi_\varepsilon$ is a UCPT map.
We will show that $\ffi_\varepsilon$ is an extremal point in the set of all UCPT maps
for sufficiently small $\varepsilon$.
Then $\ffi_\varepsilon$ goes to $\ffi$ if $\varepsilon$ goes to $0$
so that we can prove the density of 
extremal diagonal marginal tracial states of rank $a$.

To this end, we need to show that $\{w_i w_j^*\}$ is linearly independent.
Since $w_i w_j^*$ is a diagonal matrix in $M_n({\mathbb C})$,
we can consider that $w_i w_j^*$ is a vector in ${\mathbb C}^{n}$ and
let $W = [ w_i w_j^* ]_{i,j=1}^a$ be a $(n, a^2)$ matrix.
Similarly, let $U = [u_i u_j^*]_{i,j=1}^a$ and $ V=[v_i v_j^*]_{i,j=1}^a$. 
Then $\{w_i w_j^*\}$ is linearly independent if and only if
$W^* W$ is invertible in $M_{a^2}({\mathbb C})$.

Since $w_i w_j^* = (u_i + \varepsilon v_i)(u_j + \varepsilon v_j)^*
= u_i u_j^* + \varepsilon (v_i u_j^* + u_i v_j^*) + \varepsilon^2 v_i v_j^*$,
we have
\[
W = U + \varepsilon X +\varepsilon^2 V,
\]
where $X= [v_i u_j^* + u_i v_j^*]_{i,j=1}^a$, so that
\[
W^* W = U^* U + \varepsilon(U^*X + X^*U)
+ \varepsilon^2( U^*V + V^*U + X^*X)
+\varepsilon^3 (X^*V + V^* X) +\varepsilon^4 V^*V.
\]
Since $\{v_i v_j^*\}$ is linearly independent, $V^*V$ is invertible. Hence 
$W^*W$ is invertible for sufficiently small 
$\varepsilon$. Indeed,
the determinant of $W^*W$ is a polynomial of $\varepsilon$ of degree $4a^2$
and the coefficient of $\varepsilon^{4a^2}$ is the determinant of $V^*V$ 
which is not zero.
Therefore the equation $| W^*W| =0$ has at most $4a^2$ solutions and
$|W^*W|$ is not zero for sufficiently small $\varepsilon$.
\end{proof}


\section*{Acknowledgement}
The author is deeply grateful to Professor Geoffrey L. Price for helpful discussions.
The proof of Theorem \ref{thm3.3} is given by him.
The author also would like to express his gratitude to Proffessor Sh\^oichir\^o Sakai
for his useful advices.


\end{document}